\documentclass[aos,MSNbibl,nameyear,dvips]{arximspdf}

\doi{10.1214/11-AOS916}
\volume{39}
\issue{5}
\pubyear{2011}
\firstpage{2585}
\lastpage{2606}

\makeatletter

\renewcommand{\mid}{|}

\newcommand{\eqref}[1]{(\ref{#1})}

\newcommand{\norm}[1]{\Vert{#1}\Vert}

\newtheorem{proposition}{Proposition}
\newtheorem{lemma}{Lemma}
\newtheorem{theorem}{Theorem}

\newproclaim{remark}{Remark}

\newcommand{\X}{{\mathsf{X}}}
\newcommand{\Y}{{\mathsf{Y}}}

\makeatother

\begin{document}
\begin{frontmatter}

\title{A spectral analytic comparison of trace-class data
augmentation algorithms and their sandwich~variants}
\runtitle{Spectral analysis of data augmentation}

\begin{aug}
\author[A]{\fnms{Kshitij} \snm{Khare}\ead[label=e1]{kdkhare@stat.ufl.edu}}
\and
\author[A]{\fnms{James P.} \snm{Hobert}\corref{}\thanksref{t1}\ead[label=e2]{jhobert@stat.ufl.edu}}
\runauthor{K. Khare and J. P. Hobert}
\affiliation{University of Florida}
\address[A]{Department of Statistics \\
University of Florida \\
Gainesville, Florida 32611\\
USA\\
\printead{e1}\\
\hphantom{E-mail: }\printead*{e2}} %adresu isvedimo komanda gale!
\end{aug}

\thankstext{t1}{Supported by NSF Grant DMS-08-05860.}

% HISTORY:
\received{\smonth{3} \syear{2011}}
\revised{\smonth{8} \syear{2011}}

% ABSTRACT
%
\begin{abstract}
The data augmentation (DA) algorithm is a widely used Markov chain
Monte Carlo algorithm that is easy to implement but often suffers
from slow convergence. The sandwich algorithm is an alternative
that can converge much faster while requiring roughly the same
computational effort per iteration. Theoretically, the sandwich
algorithm always converges at least as fast as the corresponding DA
algorithm in the sense that $\norm{K^*} \le\norm{K}$, where $K$ and
$K^*$ are the Markov operators associated with the DA and sandwich
algorithms, respectively, and \mbox{$\Vert\cdot\Vert$} denotes operator norm.
In this paper, a substantial refinement of this operator norm
inequality is developed. In particular, under regularity conditions
implying that $K$ is a trace-class operator, it is shown that $K^*$
is also a positive, trace-class operator, and that the spectrum of
$K^*$ dominates that of $K$ in the sense that the ordered elements
of the former are all less than or equal to the corresponding
elements of the latter. Furthermore, if the sandwich algorithm is
constructed using a group action, as described by
Liu and Wu [\textit{J. Amer. Statist. Assoc.} \textbf{94} (1999)
1264--1274] and Hobert and Marchev
[\textit{Ann. Statist.} \textbf{36} (2008)
532--554], then there is strict
inequality between at least one pair of eigenvalues. These results
are applied to a new DA algorithm for Bayesian quantile regression
introduced by Kozumi and Kobayashi
[\textit{J. Stat. Comput. Simul.} \textbf{81} (2011) 1565--1578].
\end{abstract}

% KEYWORDS
%

\end{frontmatter}

%se1 #&#
\section{Introduction}
\label{secintro}

Suppose that $f_X\dvtx \X\rightarrow[0,\infty)$ is an intractable
probability density that we would like to explore. Consider a data
augmentation (DA) algorithm [\citet{tannwong1987},
\citet{liuwongkong1994}]
based on the joint density $f\dvtx \X\times\Y\rightarrow[0,\infty)$,
which of course must satisfy
\[
\int_\Y f(x,y) \nu(dy) = f_X(x) .
\]
We are assuming here that $\X$ and $\Y$ are two sets equipped with
countably generated $\sigma$-algebras, and that $f(x,y)$ is a density
with respect to $\mu\times\nu$. The Markov chain underlying the DA
algorithm, which we denote by $\{X_n\}_{n=0}^\infty$, has Markov
transition density (Mtd) given by
\[
k(x'|x) = \int_\Y f_{X|Y}(x'|y) f_{Y|X}(y|x) \nu(dy) .
\]
In other words, $k(\cdot|x)$ is the density of $X_{n+1}$, given that
$X_n=x$. It is well known and easy to see that the product $k(x'|x)
f_X(x)$ is symmetric in $(x,x')$, that is, the DA Markov chain is
reversible with respect to $f_X$. (We assume throughout that all
Markov chains on the target space, $\X$, are Harris ergodic, that is,
irreducible, aperiodic and Harris recurrent.) Of course, the DA
Markov chain can be simulated by drawing alternately from the two
conditional densities defined by $f(x,y)$. If the current state is
$X_n=x$, then~$X_{n+1}$ is simulated in two steps: draw $Y \sim
f_{Y|X}(\cdot|x)$, call the result $y$, and then draw $X_{n+1} \sim
f_{X|Y}(\cdot|y)$.

Like its cousin the EM algorithm, the DA algorithm can be very slow to
converge. A powerful method for speeding up the DA algorithm was
discovered independently by \citet{liuwu1999} and
\citet{mengvand1999}. The basic idea behind the method (called
``PX-DA'' by Liu and Wu and ``marginal augmentation'' by
Meng and van Dyk) is to introduce a (low-dimensional)
parameter into the joint density $f(x,y)$ that is not identifiable in
the target, $f_X$. This allows for the construction of an entire
class of viable DA algorithms, some of which may converge much faster
than the original. Here is a brief description of the method in the
context where $\X$ and $\Y$ are Euclidean spaces, and $f(x,y)$ is a
density with respect to the Lebesgue measure. Suppose that for each $g$
in some set $G$, there is a function \mbox{$t_g \dvtx \Y\rightarrow\Y$} that is
one-to-one and differentiable. Consider a parametric family of
densities (indexed by $g$) given by $\tilde{f}(x,y;g) = f (x,
t_g(y) ) |J_g(y)|$, where $J_g(z)$ is the Jacobian of the
transformation $z=t_g^{-1}(y)$. Note that $\int_{\Y} \tilde{f}(x,y;g)
\,dy = f_X(x)$, so $g$ is not identifiable in $f_X$. Now fix a
``working prior'' density on~$g$, call it~$r(g)$, and define a joint
density on $\X\times\Y$ as follows:
\[
f_r(x,y) = \int_G \tilde{f}(x,y;g) r(g) \,dg .
\]
Clearly, the $x$-marginal of $f_r(x,y)$ is the target, $f_X$. Thus,
each working prior leads to a new DA algorithm that is potentially
better than the original one based on $f(x,y)$. \citet{liuwu1999},
\citet{mengvand1999} and \citet{vandmeng2001} find the working
priors that lead to particularly fast algorithms.

Of course, one iteration of the DA algorithm based on $f_r(x,y)$ could
be simulated using the usual two-step method (described above) which
entails drawing from the two conditional densities defined by
$f_r(x,y)$. However, \citet{liuwu1999} showed that this simulation
can also be accomplished using a \textit{three-step} procedure in
which the first and third steps are draws from~$f_{Y|X}(\cdot|x)$ and
$f_{X|Y}(\cdot|y)$, respectively, and the middle step involves a~%
single move according to a Markov chain on the space $\Y$ that has
invariant density~$f_Y(y)$. In this paper, we study a generalization
of Liu and Wu's (\citeyear{liuwu1999}) three-step
procedure that was introduced by
\citet{hobemarc2008} and is now described.

Suppose that $R(y,dy')$ is any Markov transition function (Mtf) on
$\Y$ that is reversible with respect to $f_Y(y) \nu(dy)$, that is,
$R(y,dy') f_Y(y) \nu(dy) = R(y',dy)\times f_Y(y') \nu(dy')$.
Consider a new Mtd given by
%
%e1 #&#
%
\begin{equation}
\label{eqmtdsandwich}
k^*(x'|x) = \int_\Y\int_\Y f_{X|Y}(x'|y') R(y,dy')
f_{Y|X}(y|x) \nu(dy) .
\end{equation}
It's easy to see that $k^*(x'|x) f_X(x)$ is symmetric in $(x,x')$,
so the Markov chain defined by $k^*$, which we denote by $\{ X^*_n
\}_{n=0}^\infty$, is reversible with respect to $f_X$. If the current
state of the new chain is $X^*_n=x$, then $X^*_{n+1}$ can be simulated
using the following three-steps, which are suggested by the form of
$k^*$. Draw $Y \sim f_{Y|X}(\cdot|x)$, call the result $y$, then draw
$Y' \sim R(y,\cdot)$, call the result $y'$, and finally draw
$X^*_{n+1} \sim f_{X|Y}(\cdot|y')$. Again, the first and third steps
are exactly the two steps used to simulate the original DA algorithm.
Because the draw from $R(y,\cdot)$ is ``sandwiched'' between the draws
from the two conditional densities, \citet{yumeng2011} call this new
algorithm the ``sandwich algorithm'' and we will follow their lead.
The PX-DA/marginal augmentation method can be viewed as one particular
recipe for constructing~$R(y,dy')$. Another general method for
building $R$ is described in Section~\ref{secgroup}.

It is often possible to construct a sandwich algorithm that converges
much faster than the underlying DA algorithm while requiring roughly
the same computational effort per iteration. Examples can be found in
\citet{liuwu1999}, \citet{mengvand1999}, \citet{vandmeng2001},
\citet{marchobe2004} and \citet{hoberoyrobe2011}. What makes
this ``free lunch'' possible is the somewhat surprising fact that a
low-dimensional perturbation on the $\Y$ space can lead to a major
improvement in mixing. In fact, the chain driven by $R$ is typically
reducible, living in a small subspace of $\Y$ that is determined by
its starting value. Drawing from such an $R$ is usually much less
expensive computationally than drawing from $f_{Y|X}(\cdot|x)$ and
$f_{X|Y}(\cdot|y)$.

Empirical studies pointing to the superiority of the sandwich
algorithm abound. Unfortunately, the development of confirmatory
theoretical results has been slow. It is known that the sandwich
chain always converges at least as fast as the DA chain in the
operator norm sense. Indeed, \citet{hoberoma2011} show that
Yu and Meng's (\citeyear{yumeng2011}) Theorem 1 can
be used to show that
%
%e2 #&#
%
\begin{equation}
\label{eqnorms}
\norm{K^*} \le\norm{R} \norm{K} ,
\end{equation}
where $K$, $K^*$ and $R$ denote the usual Markov operators defined by
$k$, $k^*$ and~$R(y,dy')$, respectively, and \mbox{$\Vert\cdot\Vert$} denotes
the operator norm. (See Section~\ref{secreview} for more details as well
as references.) Of course, we would like to be able to say that
$\norm{K^*}$ is strictly smaller than $\norm{K}$, and this is
certainly the case when $\norm{R}<1$. However, the $R$'s used in
practice typically have norm 1 (because the corresponding chains are
reducible). In fact, in most applications, $R$ is reducible and
idempotent, that is, $\int_\Y R(y,dy'') R(y'',dy') = R(y,dy')$.

\citet{hoberoyrobe2011} provided a refinement of \eqref{eqnorms}
for the case in which $\Y$ is finite and $R$ is reducible and
idempotent. These authors showed that, in this case, $K$ and $K^*$
both have pure eigenvalue spectra that are subsets of $[0,1)$, and
that at most $m-1$ of the eigenvalues are nonzero, where
$|\Y|=m<\infty$. They also showed that the spectrum of $K^*$
dominates that of $K$ in the sense that $0 \le\lambda^*_i \le
\lambda_i < 1$ for all $i$, where $\lambda_i$ and $\lambda^*_i$ denote
the $i$th largest eigenvalues of $K$ and $K^*$, respectively. Note
that taking $i=1$ yields $\norm{K^*} = \lambda^*_1 \le\lambda_1 =
\norm{K}$.

In this paper we develop results that hold in the far more common
situation where $|\Y|=\infty$. First, we generalize
Hobert, Roy and Robert's (\citeyear{hoberoyrobe2011})
result by showing that the assumption that
$\Y$ is finite can be replaced by the substantially weaker assumption
that $\int_\X k(x|x) \mu(dx) < \infty$. In this more general case,
$K$ and $K^*$ still have pure eigenvalue spectra that are subsets of
$[0,1)$ and an analogous domination holds, but the number of nonzero
eigenvalues is no longer necessarily finite. Second, we show that
if $R$ is constructed using a group action, as described by
\citet{liuwu1999} and \citet{hobemarc2008}, then the domination is
strict in the sense that there exists at least one $i$ such that $0
\le\lambda^*_i < \lambda_i < 1$. Finally, we apply our results to a
new DA algorithm for Bayesian quantile regression that was recently
introduced by \citet{kozukoba2011}.

The remainder of this paper is organized as follows.
Section \ref{secreview} contains a~brief review of the relationship
between the spectral properties of Markov operators and the
convergence properties of the corresponding Markov chains. The DA and
sandwich algorithms are formally defined and compared in
Section \ref{secalgorithms}. The construction of $R$ using group
actions is discussed in Section \ref{secgroup}, and our analysis of
Kozumi and Kobayashi's DA algorithm is described in
Section \ref{secexample}.

%se2 #&#
\section{Brief review of self-adjoint Markov operators}
\label{secreview}

Let $P(x,dx')$ be a~gene\-ric Mtf on $\X$ that is reversible with
respect to $f_X(x) \mu(dx)$. Denote the Markov chain driven by $P$ as
$\Phi= \{\Phi_n\}_{n=0}^\infty$. (Note that $\Phi$ is not
necessarily a DA Markov chain.) The convergence properties of $\Phi$
can be expressed in terms of a related operator that is now defined.
Let $L^2_0(f_X)$ be the space of real-valued functions with domain
$\X$ that are square integrable and have mean zero with respect to
$f_X$. In other words, $g \in L^2_0(f_X)$ if $\int_\X g^p(x)
f_X(x) \mu(dx)$ is finite when $p=2$, and vanishes when $p=1$.
This is a Hilbert space where the inner product of $g,h \in L^2_0(f_X)$ is
defined as
\[
\langle g,h \rangle= \int_\X g(x) h(x) f_X(x) \mu(dx) ,
\]
and the corresponding norm is, of course, given by $\norm{g} = \langle
g,g \rangle^{1/2}$. Let $P\dvtx\break L^2_0(f_X) \rightarrow L^2_0(f_X)$ denote
the operator that maps $g \in L^2_0(f_X)$ to
\[
(Pg)(x) = \int_{\X} g(x') P(x,dx') .
\]
Note that $(Pg)(x)$ is simply the conditional expectation of
$g(\Phi_{n+1})$ given that $\Phi_n=x$. Reversibility of the Mtf
$P(x,dx')$ is equivalent to the operator~$P$ being self-adjoint. The
(operator) norm of $P$ is defined as
\[
\norm{P} = {\sup_{g \in L^2_{0,1}(f_X)}} \norm{Pg} ,
\]
where $L^2_{0,1}(f_X)$ in the subset of $L^2_0(f_X)$ that contains the
functions $g$ satisfying $\int_\X g^2(x) f_X(x) \mu(dx) = 1$.
It's easy to see that $\norm{P} \in[0,1]$. \citet{roberose1997}
show that $\norm{P} < 1$ if and only if $\Phi$ is geometrically
ergodic. Moreover, in the geometrically ergodic case, $\norm{P}$ can
be viewed as the asymptotic rate of convergence of $\Phi$
[see, e.g., \citet{rose2003}, page~170].

If $P$ satisfies additional regularity conditions, much more can be
said about the convergence of the corresponding Markov chain. Assume
that the operator $P$ is compact and positive, and let $\{\alpha_i\}$
denote its eigenvalues, all of which reside in $[0,1)$. The number of
eigenvalues may be finite or countably infinite (depending on the
cardinality of $\X$), but in either case there is a largest one and it
is equal to $\norm{P}$. [For a basic introduction to the spectral
properties of Markov operators, see \citet{hoberoyrobe2011}.] If
$\operatorname{tr}(P) := \sum_i \alpha_i < \infty$, then $P$ is a
\textit{trace-class} operator [see, e.g., \citet{conw1990}, page 267].
As explained in \citet{diackharsalo2008}, when $P$ is positive and
trace-class, the chain's $\chi^2$-distance to stationarity can be
written explicitly as
%
%e3 #&#
%
\begin{equation}
\label{eqrep}
\int_\X\frac{| p^n(x'|x) - f_X(x') |^2}{f_X(x')} \mu(dx')
= \sum_i \alpha_i^{2n} e_i^2(x) ,
\end{equation}
where $p^n(\cdot|x)$ denotes the density of $\Phi_n$ given that
$\Phi_0=x$, and $\{e_i\}$ is an orthonormal basis of
eigen-functions\vadjust{\goodbreak}
corresponding to $\{\alpha_i\}$. Of course, the $\chi^2$-distance
serves as an upper bound on the total variation distance. Assume that
the eigenvalues are ordered so that $\alpha_i \ge\alpha_{i+1}$, and
let $i^* = \max\{ i \in\mathbb{N}\dvtx \alpha_i = \alpha_1
\}$. Asymptotically, the term $\alpha_1^{2n} ( e_1^2(x) + \cdots
+ e_{i^*}^2(x) )$ will dominate the sum on the right-hand side of
\eqref{eqrep}. Hence, in this context it is certainly reasonable to
call $\norm{P} = \alpha_1$ the ``asymptotic rate of convergence.''
Our focus in this paper will be on DA algorithms whose Markov
operators are trace-class.

%se3 #&#
\section{Spectral comparison of the DA and sandwich algorithms}
\label{secalgorithms}

As in Section~\ref{secintro}, let $K\dvtx L_0^2(f_X) \rightarrow
L_0^2(f_X)$, $K^*\dvtx L_0^2(f_X) \rightarrow L_0^2(f_X)$ and $R\dvtx
L_0^2(f_Y) \rightarrow L_0^2(f_Y)$ denote the (self-adjoint) Markov
operators defined by the DA chain, the sandwich chain and $R(y,dy')$,
respectively. We will exploit the fact that~$K$ and $K^*$ can be
represented as products of simpler operators. Indeed, let $P_X\dvtx
L^2_0(f_Y) \rightarrow L^2_0(f_X)$ map $h \in L^2_0(f_Y)$ to
\[
(P_X h)(x) = \int_{\Y} h(y) f_{Y|X}(y|x) \nu(dy)
\]
and, analogously, let $P_Y\dvtx L^2_0(f_X) \rightarrow L^2_0(f_Y)$ map $g
\in L^2_0(f_X)$ to
\[
(P_Y g)(y) = \int_{\X} g(x) f_{X|Y}(x|y) \mu(dx) .
\]
It is easy to see that $K = P_X P_Y$ and $K^* = P_X R P_Y$. This
representation of~$K$ was used in \citet{diackharsalo2008}.

Again, as in Section \ref{secintro}, let $f\dvtx \X\times\Y\rightarrow
[0,\infty)$ be the joint density that defines the DA Markov chain.
Throughout the next two sections, we assume that $f$ satisfies the
following condition:
%
%e4 #&#
%
\begin{equation}
\label{eqct}
\int_\X\int_\Y f_{X|Y}(x|y) f_{Y|X}(y|x) \nu(dy) \mu(dx) <
\infty.
\end{equation}
\citet{buja1990} shows that \eqref{eqct} implies that $K$ is a
trace-class operator. It is clear that~\eqref{eqct} holds if $\X$
and/or $\Y$ has a finite number of elements. However,~\eqref{eqct}
can also hold in situations where $|\X| = |\Y| = \infty$. Indeed, in
Section \ref{secexample} we establish that \eqref{eqct} holds for a
DA algorithm for Bayesian quantile regression where $\X$ and $\Y$ are
both uncountable. On the other hand, \eqref{eqct} certainly does not
hold for all DA algorithms. For example, \eqref{eqct} cannot hold if
the DA chain is not geometrically ergodic (because subgeometric chains
cannot be trace-class). Simple examples of subgeometric DA chains can
be found in \citet{paparobe2008} and \citet{tan2008}, Chapter 4.

Condition \eqref{eqct} has appeared in the Markov chain Monte Carlo
literature before. It is exactly the bivariate version of
Liu, Wong and Kong's (\citeyear{liuwongkong1995})
``Condition (b)'' and it was also employed
by \citet{schecarl1992}. Unfortunately, there does not appear to be
any simple, intuitive interpretation of \eqref{eqct} in terms of the
joint density $f(x,y)$ or the corresponding Markov chain. In fact,
referring to their Condition (b), \citet{liuwongkong1995} state that
``It is standard but not easy to check and understand.''

Our analysis of the DA and sandwich algorithms rests heavily upon a~%
useful \textit{singular value decomposition} of $f(x,y)$ whose
existence is implied by~\eqref{eqct}. In particular,
\citet{buja1990} shows that if \eqref{eqct} holds, then
%
%e5 #&#
%
\begin{equation}
\label{eqsvd}
\frac{f(x,y)}{f_X(x) f_Y(y)} = \sum_{i=0}^\infty\beta_i
g_i(x) h_i(y) ,
\end{equation}
where:
\begin{itemize}
\item$\beta_0 = 1$, $g_0 \equiv1$, $h_0 \equiv1$.
\item$\{g_i\}_{i=0}^\infty$ and $\{h_i\}_{i=0}^\infty$ form
orthonormal bases of $L^2(f_X)$ and $L^2(f_Y)$, respectively.
\item$\beta_i \in[0,1]$, and $\beta_i \le\beta_{i-1}$ for all
$i \in\mathbb{N}$.
\item$\int_\X\int_\Y g_i(x) h_j(y) f(x,y) \nu(dy)
\mu(dx) = 0$ if $i \ne j$.
\end{itemize}

A few remarks about notation are in order. First, we state all
results for the case $|\X| = |\Y| = \infty$, and leave it to the
reader to make the obvious, minor modifications that are required when
one or both of the spaces are finite. For example, in the singular
value decomposition above, if one or both of the spaces are finite,
then one or both of the orthonormal bases would have only a finite
number of elements, etc. Second, we will let $\langle\cdot,\cdot
\rangle$ and \mbox{$\Vert\cdot\Vert$} do double duty as inner product and norm
on both $L^2_0(f_X)$ and $L^2_0(f_Y)$. However, the norms of
operators whose domains and ranges differ will be subscripted. The
following result can be gleaned from calculations in
\citet{buja1990}, but we present a proof in Appendix \ref{applemma1}
for completeness.
%
%le1 #&#
%
\begin{lemma}
\label{lemlambda1}
Assume that \eqref{eqct} holds and let $\lambda_1 \ge\lambda_2 \ge
\lambda_3 \ge\cdots$ denote the eigenvalues of $K$, which reside in
the set $[0,1)$. For each $i \in\mathbb{N}$, $P_X h_i = \beta_i
g_i$ and $P_Y g_i = \beta_i h_i$. Moreover,
\[
\norm{P_X}_{L^2_0(f_Y) \rightarrow L^2_0(f_X)} =
\norm{P_Y}_{L^2_0(f_X) \rightarrow L^2_0(f_Y)} = \beta_1
\]
and $\lambda_i = \beta_i^2$.
\end{lemma}

Here is the first of our two main results.
%
%th1 #&#
%
\begin{theorem}
\label{thmcompare}
Assume that \eqref{eqct} holds and that $R$ is idempotent with
$\norm{R}=1$. Define $l = \max\{ i \in\mathbb{N}\dvtx \beta_i =
\beta_1 \}$ and $N = \{ i \in\mathbb{N}\dvtx \beta_i > 0
\}$. Then:
{\renewcommand\thelonglist{\arabic{longlist}}
\begin{longlist}
\item$K^*$ is a positive, trace-class operator.
\item$\lambda^*_i \le\lambda_i$ for all $i \in\mathbb{N}$, where
$\{\lambda_i\}_{i=1}^\infty$ and $\{\lambda^*_i\}_{i=1}^\infty$
denote the (ordered) eigenvalues of $K$ and $K^*$, respectively.
\item$\lambda^*_i = \lambda_i$ for all $i \in\mathbb{N}$ if and only
if $Rh_i = h_i$ for every $i \in N$.
\item A necessary and sufficient condition for $\norm{K^*} < \norm{K}$
is that the only $a = ( a_1,\ldots,a_l) \in\mathbb{R}^l$ for which
%
%e6 #&#
%
\begin{equation}
\label{eqkey}
R \sum_{i=1}^l a_i h_i = \sum_{i=1}^l a_i h_i
\end{equation}
is the zero vector in $\mathbb{R}^l$.
\end{longlist}}
\end{theorem}
%
%re1 #&#
%
\begin{remark}
Part (3) can be rephrased as follows: $\operatorname{tr}(K^*) = \operatorname{tr}(K)$
if and only if $Rh_i = h_i$ for every $i \in N$. Also, note that
$\norm{K^*} = \lambda^*_1$ and $\norm{K} = \lambda_1$.
\end{remark}
\begin{pf*}{Proof of Theorem \ref{thmcompare}}
We begin by noting that for $g \in L^2_0(f_X)$ and $h \in
L^2_0(f_Y)$, we have $\langle P_X h,g \rangle= \langle h, P_Yg
\rangle$. Hence,
\[
\langle K^* g,g \rangle= \langle P_X R P_Y g,g \rangle= \langle R
P_Y g, P_Y g \rangle= \langle R^{1/2} P_Y g, R^{1/2}
P_Y g \rangle\ge0 ,
\]
which shows that $K^*$ is positive. Since $K$ is trace-class, it
follows from Lem\-ma~\ref{lemlambda1} that $\operatorname{tr}(K) =
\sum_{i=1}^\infty\beta_i^2 < \infty$. Now, since
$\{g_i\}_{i=1}^\infty$ is an orthonormal basis for $L^2_0(f_X)$, we
have
\begin{eqnarray*}
\operatorname{tr}(K^*) &=& \sum_{i=1}^\infty\langle K^* g_i, g_i
\rangle= \sum_{i=1}^\infty\langle P_X R  P_Y g_i, g_i
\rangle= \sum_{i=1}^\infty\langle R P_Y g_i, P_Y g_i
\rangle\\
&=& \sum_{i=1}^\infty\beta_i^2 \langle R h_i, h_i
\rangle\le\sum_{i=1}^\infty\beta_i^2 = \operatorname{tr}(K) ,
\end{eqnarray*}
where the inequality follows from the fact that $\langle R h_i,
h_i \rangle\le1$. Thus, $K^*$ is trace-class. Moreover, it is
clear that $\operatorname{tr}(K^*) = \operatorname{tr}(K)$ if and only if $
\langle R h_i, h_i \rangle= 1$ whenever $\beta_i > 0$. Since
$R$ is idempotent with norm 1, it is a projection
[\citet{conw1990}, page 37]. Thus, for any $h \in L^2_0(f_Y)$,
$\langle
Rh,h \rangle= \langle h,h \rangle\Rightarrow Rh = h$. [Indeed,
$\langle h,h \rangle= \langle Rh,h \rangle+ \langle(I-R)h,h
\rangle$, so $\langle Rh,h \rangle= \langle h,h \rangle\Rightarrow
\langle(I-R)h,h \rangle= \langle(I-R)h,(I-R)h \rangle= 0$.]
Consequently, $\operatorname{tr}(K^*) = \operatorname{tr}(K)$ if and only if $Rh_i =
h_i$ for every $i$ such that $\beta_i > 0$. This takes care of (3).

Now, note that $K-K^* = P_X(I-R)P_Y$ is positive since
\[
\langle P_X(I-R)P_Yg,g \rangle= \langle(I-R)P_Yg,P_Yg \rangle=
\langle(I-R)P_Yg,(I-R)P_Yg \rangle\ge0.
\]
Therefore, for any nonnull $g \in L^2_0(f_X)$, we have
\[
\frac{\langle K^* g,g \rangle}{\langle g,g \rangle} \le\frac
{\langle
Kg,g \rangle}{\langle g,g \rangle}.
\]
Now, for any $i \in\mathbb{N}$, the Courant--Fischer--Weyl minmax
characterization of eigenvalues of compact, positive, self-adjoint
operators [see, e.g., \citet{voss2003}] yields
\[
\lambda^*_i = \min_{\operatorname{dim}(V) = i-1} \max_{g \in V^\perp, g
\ne0} \frac{\langle K^* g,g \rangle}{\langle g,g \rangle} \le
\min_{\operatorname{dim}(V) = i-1} \max_{g \in V^\perp, g \ne0}
\frac{\langle Kg,g \rangle}{\langle g,g \rangle} = \lambda_i ,
\]
where $V$ denotes a subspace of $L^2_0(f_X)$, and $\operatorname{dim}(V)$ is
its dimension. This proves~(2).

All that remains is (4). Assume there exists a nonzero $a$ such that
\eqref{eqkey} holds. We will show that $\norm{K^*} = \norm{K}$.
Since we know that $\norm{K^*} \le\norm{K} = \beta_1^2$, it suffices
to identify a function $g \in L_0^2(f_X)$ such that $\norm{K^* g} =
\beta_1^2 \norm{g}$. If we take $g = a_1 g_1 + \cdots+ a_l g_l$,
then
\[
\norm{K^* g} = \Biggl\| K^* \sum_{i=1}^l a_i g_i \Biggr\| =
\Biggl\| P_X R P_Y \sum_{i=1}^l a_i g_i \Biggr\| = \Biggl\| P_X R
\sum_{i=1}^l a_i \beta_i h_i \Biggr\| .
\]
But $\beta_1 = \cdots= \beta_l$, and, hence,
\[
\norm{K^* g} = \beta_1 \Biggl\| P_X R \sum_{i=1}^l a_i h_i
\Biggr\|
= \beta_1 \Biggl\| P_X \sum_{i=1}^l a_i h_i \Biggr\| = \beta_1^2
\Biggl\| \sum_{i=1}^l a_i g_i \Biggr\| = \beta^2_1 \norm{g} .
\]
The second half of the proof is by contradiction. Assume that the
only $a \in\mathbb{R}^l$ for which \eqref{eqkey} holds is the zero
vector, and assume also that $\norm{K^*} = \norm{K} = \beta_1^2$. By
completeness of the Hilbert space, $L^2_0(f_X)$, there exists a
nontrivial function $g \in L^2_0(f_X)$ such that $\norm{K^* g} =
\beta_1^2 \norm{g}$. The rest of the argument differs depending upon
whether $g$ is in the span of $\{ g_1,\ldots,g_l \}$ or not.

\textit{Case} I: Assume that $g = \sum_{i=1}^l a_i g_i$ for
some nonzero $a \in\mathbb{R}^l$. Using the results above, we have
\[
\norm{K^* g} = \| P_X R P_Y g \| \le\beta_1 \| R
P_Y g \| = \beta_1^2 \Biggl\| R \sum_{i=1}^l a_i h_i \Biggr\|.
\]
But $R$ is a projection, so $Rh \ne h \Rightarrow\norm{Rh} \ne
\norm{h}$. Hence, $R \sum_{i=1}^l a_i h_i \ne\sum_{i=1}^l a_i h_i$
in conjunction with $\norm{R}=1$ yields
\[
\Biggl\| R \sum_{i=1}^l a_i h_i \Biggr\| < \Biggl\| \sum_{i=1}^l a_i
h_i \Biggr\| = \sqrt{ \sum_{i=1}^l a_i^2} = \norm{g}.
\]
Thus, $\norm{K^* g} < \beta_1^2 \norm{g}$, which is a
contradiction.

\textit{Case} II: Assume that $g$ is not in the span of $\{
g_1,\ldots,g_l \}$. In other words, $g = \sum_{i=1}^\infty b_i g_i$
where at least one term in the sequence $\{ b_{l+1},b_{l+2},\ldots\}$
is nonzero. Then,
\[
\norm{P_Y g} = \Biggl\| P_Y \sum_{i=1}^\infty b_i g_i \Biggr\| =
\Biggl\| \sum_{i=1}^\infty b_i \beta_i h_i \Biggr\| =
\sqrt{\sum_{i=1}^\infty b_i^2 \beta_i^2} < \sqrt{\beta_1^2
\sum_{i=1}^\infty b_i^2} = \beta_1 \norm{g} .
\]
It follows that
\[
\norm{K^* g} \le\norm{P_X}_{L^2_0(f_Y) \rightarrow L^2_0(f_X)}
\norm{R} \norm{P_Y g} < \beta_1^2 \norm{g} ,
\]
and, again, this is a contradiction.
\end{pf*}

%se4 #&#
\section{Using a group action to construct $R$}
\label{secgroup}

Following \citet{liuwu1999} and \citet{liusaba2000},
\citet{hobemarc2008} introduced and studied a general method for
constructing practically useful versions of~$R(y,\allowbreak dy')$ using group
actions. For the remainder of this section, assume that $\X$ and~$\Y$
are locally compact, separable metric spaces equipped with their Borel
$\sigma$-algebras. Suppose that $G$ is a third locally compact,
separable metric space that is also a topological group. As usual,
let $e$ denote the identity element of the group. Also, let
$\mathbb{R}_+ = (0,\infty)$. Any continuous function $\chi\dvtx G
\rightarrow\mathbb{R}_+$ such that $\chi(g_1 g_2)$ = $\chi(g_1)
\chi(g_2)$ for all $g_1,g_2 \in G$ is called a \textit{multiplier}
[\citet{eato1989}]. Clearly, a multiplier must satisfy $\chi(e) = 1$
and $\chi( g^{-1} ) = 1/\chi(g)$. One important multiplier
is the \textit{modular function}, $\Delta$, which relates the
left-Haar and right-Haar measures on $G$. Indeed, if we denote these
measures by $\omega_l(\cdot)$ and $\omega_r(\cdot)$, then
$\omega_r(dg) = \Delta( g^{-1} ) \omega_l(dg)$. Groups for
which $\Delta\equiv1$ are called unimodular groups.

An example (that will be used later in Section \ref{secexample}) is
the multiplicative group, $\mathbb{R}_+$, where group\vspace*{1pt} composition is
defined as multiplication, the identity element is $e=1$ and $g^{-1}
= 1/g$. This group is unimodular with Haar measure given by
$\omega(dg) = dg/g$ where $dg$ denotes the Lebesgue measure on
$\mathbb{R}_+$.

Let $F\dvtx G \times\Y\rightarrow\Y$ be a continuous function
satisfying $F(e,y)=y$ and $F(g_1g_2,y) = F(g_1,F(g_2,y))$ for all
$g_1,g_2 \in G$ and all $y \in\Y$. The function~$F$ represents $G$
acting topologically on the left of $\Y$ and, as is typical, we
abbreviate $F(g,y)$ with $gy$. Now suppose there exists a multiplier,
$\chi$, such~that\vspace*{-1pt}
\[
\chi(g) \int_\Y\phi(gy) \nu(dy) = \int_\Y\phi(y) \nu(dy)\vspace*{-1pt}
\]
for all $g \in G$ and all integrable $\phi\dvtx \Y\rightarrow
\mathbb{R}$. Then the measure $\nu$ is called \textit{relatively}
(\textit{left}) \textit{invariant} with multiplier $\chi$. For example, suppose that
\mbox{$\Y= \mathbb{R}^m$}, $\nu(dy)$ is the Lebesgue measure, $G$ is the
multiplicative group described above, and the group action is defined
to be scalar multiplication, that is, $gy = (gy_1,gy_2,\ldots,gy_m)$.
Then $\nu(dy)$ is relatively invariant with multiplier $\chi(g) =
g^m$. Indeed,\vspace*{-1pt}
\[
g^m \int_{\mathbb{R}^m} \phi(gy) \nu(dy) = \int_{\mathbb{R}^m}
\phi(y) \nu(dy) .\vspace*{-1pt}
\]

We now explain how the group action is used to construct $R$. Define\vspace*{-1pt}
\[
m(y) = \int_G f_Y(gy) \chi(g) \omega_l(dg) .\vspace*{-1pt}
\]
Assume that $m(y)$ is positive for all $y \in\Y$ and finite for
$\nu$-almost all $y \in\Y$. For the remainder of this section, we
assume that $R\dvtx L^2_0(f_Y) \rightarrow L^2_0(f_Y)$ is the operator
that maps $h(y)$ to\vspace*{-1pt}
\[
(Rh)(y) = \frac{1}{m(y)} \int_G h(gy) f_Y(gy) \chi(g)
\omega_l(dg) .\vspace*{-1pt}
\]
\citet{hobemarc2008} show that $R$ is a self-adjoint, idempotent
Markov operator on $L^2_0(f_Y)$. The corresponding\vadjust{\goodbreak} Markov chain on
$\Y$ evolves as follows. If the current state is $y$, then the
distribution of the next state is that of $gy$, where $g$ is a random
element from $G$ whose density is
%
%e7 #&#
%
\begin{equation}
\label{eqgroupden}
\frac{f_Y(gy)\chi(g)}{m(y)} \omega_l(dg) .
\end{equation}
Therefore,\vspace*{1pt} we can move from $X^*_n=x$ to $X^*_{n+1}$ as follows: draw
$Y \sim f_{Y|X}(\cdot|x)$, call the result $y$, then draw $g$ from the
density \eqref{eqgroupden} and set $y'=gy$, and finally draw
$X^*_{n+1} \sim f_{X|Y}(\cdot|y')$.\vspace*{1pt}

\citet{hobemarc2008} also show that, if $\{Y_n\}_{n=0}^\infty$
denotes the Markov chain defined by $R$, then conditional on $Y_0=y$,
$\{Y_n\}_{n=1}^\infty$ are i.i.d. Thus, either $\{Y_n\}_{n=1}^\infty$
are i.i.d. from $f_Y$, or the chain is reducible.
%
%le2 #&#
%
\begin{lemma}
\label{leminvariantproj}
If $\norm{R}=1$, then the Markov operator $R$ is a projection onto
the space of functions that are invariant under the group action,
that is, $h$ is in the range of $R$ if and only if $h(gy) = h(y)$
for all $g \in G$ and all $y \in\Y$.
\end{lemma}
\begin{pf}
$\!\!\!$First, assume that $h(gy) = h(y)$ for all $g \in G$ and
all $y \in\Y$.~Then
\begin{eqnarray*}
(Rh)(y) & = & \frac{1}{m(y)} \int_G h(gy) f_Y(gy) \chi(g)
\omega_l(dg) \\
& = &\frac{h(y)}{m(y)} \int_G f_Y(gy) \chi(g)
\omega_l(dg) = h(y) .
\end{eqnarray*}
To prove the necessity, we require two results that were used
repeatedly by \citet{hobemarc2008}. First,
%
%e8 #&#
%
\begin{equation}
\label{eqh-m2}
\chi(g) m(gy) = \Delta( g^{-1} ) m(y) .
\end{equation}
Second, if $\tilde{g} \in G$ and $\phi\dvtx G \rightarrow\mathbb{R}$ is
integrable with respect to $\omega_l$, then
%
%e9 #&#
%
\begin{equation}
\label{eqh-m1}
\int_G \phi(g \tilde{g}^{-1} ) \omega_l(dg) = \Delta(
\tilde{g} ) \int_G \phi(g) \omega_l(dg) .
\end{equation}
Now, fix $h \in L_0^2(f_Y)$ and $g' \in G$, and note that
\begin{eqnarray*}
(Rh)(g'y) & = & \frac{1}{m(g'y)} \int_G h(gg'y) f_Y(gg'y)
\chi(g) \omega_l(dg) \\
& = & \frac{1}{\chi(g') m(g'y)} \int_G
h(gg'y) f_Y(gg'y) \chi(gg') \omega_l(dg) \\
& = & \frac{\Delta( g'^{-1} )}{\chi(g') m(g'y)} \int_G h(gy)
f_Y(gy) \chi(g) \omega_l(dg) \\
& = & \frac{1}{m(y)} \int_G
h(gy) f_Y(gy) \chi(g) \omega_l(dg) \\
& = & (Rh)(y) ,
\end{eqnarray*}
where the third and fourth equalities are due to \eqref{eqh-m1} and
\eqref{eqh-m2}, respectively.\vadjust{\goodbreak}~%
\end{pf}

We now show that when $R$ is constructed using the group action recipe
described above, there is at least one eigenvalue of $K^*$ that is
strictly smaller than the corresponding eigenvalue of $K$. To get a
strict inequality, we must rule out trivial cases in which the DA and
sandwich algorithms are the same. For example, if we take $G$ to be
the subgroup of the multiplicative group that contains only the point
$\{1\}$, then element-wise multiplication of $y \in\mathbb{R}^m$ by
$g$ has no effect and the sandwich algorithm is the same as the DA
algorithm. More generally, if
%
%e10 #&#
%
\begin{equation}
\label{eqcanthold}
f_{X|Y}(x|y) = f_{X|Y}(x|gy) \qquad\forall g \in G, x \in\X, y
\in\Y,
\end{equation}
then the Mtd of the sandwich chain can be expressed as
\begin{eqnarray*}
k^*(x'|x) & = & \int_\Y\int_G f_{X|Y}(x'|gy) \biggl[
\frac{f_Y(gy)\chi(g)}{m(y)} \omega_l(dg) \biggr] f_{Y|X}(y|x)
\nu(dy) \\
& = & \int_\Y\int_G f_{X|Y}(x'|y) \biggl[
\frac{f_Y(gy)\chi(g)}{m(y)} \omega_l(dg) \biggr] f_{Y|X}(y|x)
\nu(dy) \\ & = & \int_\Y f_{X|Y}(x'|y) f_{Y|X}(y|x) \nu(dy) \\
& = & k(x'|x) .
\end{eqnarray*}
Thus, \eqref{eqcanthold} implies that the DA and sandwich algorithms
are \textit{exactly the same} and, consequently, $\operatorname{tr}(K^*) =
\operatorname{tr}(K)$. In fact, as the next result shows,~%
\eqref{eqcanthold} is also necessary for $\operatorname{tr}(K^*) =
\operatorname{tr}(K)$.
%
%th2 #&#
%
\begin{theorem}
\label{thmatleastone}
If \eqref{eqcanthold} does not hold, then $\operatorname{tr}(K^*) <
\operatorname{tr}(K)$, so at least one eigenvalue of $K^*$ is strictly
smaller than the corresponding eigenvalue of $K$.
\end{theorem}
\begin{pf}
It is enough to show that $\operatorname{tr}(K^*) = \operatorname{tr}(K)$ implies
\eqref{eqcanthold}. Recall that $N = \{ i \in\mathbb{N} \dvtx
\beta_i > 0 \}$. By Theorem \ref{thmcompare}, $\operatorname{tr}(K^*)
= \operatorname{tr}(K)$ implies that $Rh_i = h_i$ for every $i \in N$. By
Lemma \ref{leminvariantproj}, if $Rh_i = h_i$ for every $i \in N$,
then every member of the set $\{ h_i\dvtx i \in N \}$ is
invariant under the group action. Now, using the singular value
decomposition, we see that for every $g \in G, x \in\X, y \in\Y$,
we have
\begin{eqnarray*}
f_{X|Y}(x|y) &=& \sum_{j=0}^\infty\beta_j g_j(x) h_j(y) f_X(x) \\
&=&
\sum_{j=0}^\infty\beta_j g_j(x) h_j(gy) f_X(x) \\
&=& f_{X|Y}(x|gy) .
\end{eqnarray*}
\upqed\end{pf}

In practice, $f_{X|Y}(x|y)$ is known exactly and it's easy to verify
that \eqref{eqcanthold} does not hold. An example is given in the
next section.\vadjust{\goodbreak}

It is important to note that, while Theorem \ref{thmatleastone}
guarantees strict inequality between at least one pair of eigenvalues
of $K$ and $K^*$, it does not preclude equality of $\lambda_1$ and
$\lambda^*_1$. Thus, we could still have $\norm{K} = \norm{K^*}$. We
actually believe that one would have to be quite unlucky to end up in
a situation where $\norm{K} = \norm{K^*}$. To keep things simple,
suppose that the largest eigenvalue of $K$ is unique. According to
Theorem \ref{thmcompare} and Lemma \ref{leminvariantproj},
$\norm{K} = \norm{K^*}$ if and only if $h_1$ [from \eqref{eqsvd}] is
invariant under the group action. This seems rather unlikely given
that the choice of group action is usually based on simplicity and
convenience. This is borne out in the toy examples analyzed by
\citet{hoberoyrobe2011} where there is strict inequality among
\textit{all} pairs of eigenvalues.

Recall from Section \ref{secintro} that the PX-DA/marginal
augmentation algorithm is based on a class of transformations $t_g
\dvtx
\Y\rightarrow\Y$, for $g \in G$. This class can sometimes be used
[as the function $F(g,y)$] to construct an $R$ as described above, and
when this is the case, the resulting sandwich algorithm is the same as
the optimal limiting PX-DA/marginal augmentation algorithm
[\citet{liuwu1999}, \citet{mengvand1999},
\citet{hobemarc2008}].

%se5 #&#
\section{A DA algorithm for Bayesian quantile regression}
\label{secexample}

Suppose $Z_1,Z_2,\ldots,\allowbreak Z_m$ are independent random variables such that
$Z_i = x_i^T \beta+ \varepsilon_i$ where $x_i \in\mathbb{R}^p$ is a
vector of known covariates associated with $Z_i$, $\beta\in
\mathbb{R}^p$ is a vector of unknown regression coefficients, and
$\varepsilon_1,\ldots,\varepsilon_m$ are i.i.d. errors with common density given
by
\[
d(\varepsilon;r) = r(1-r) \bigl[ e^{(1-r)\varepsilon}
I_{\mathbb{R}_-}(\varepsilon) + e^{-r \varepsilon}
I_{\mathbb{R}_+}(\varepsilon) \bigr] ,
\]
where $r \in(0,1)$. This error density, called the asymmetric
Laplace density, has $r$th quantile equal to zero. Note that when
$r=1/2$, it is the usual Laplace density with location and scale equal
to 0 and $1/2$, respectively.

If we put a flat prior on $\beta$, then the product of the likelihood
function and the prior is equal to $r^m(1-r)^m s(\beta,z)$, where
\[
s(\beta,z) := \prod_{i=1}^m \bigl[ e^{(1-r) ( z_i - x_i^T \beta
)} I_{\mathbb{R}_-}(z_i - x_i^T \beta) + e^{-r ( z_i - x_i^T
\beta)} I_{\mathbb{R}_+}(z_i - x_i^T \beta) \bigr] .
\]
If $s(\beta,z)$ is normalizable, that is, if
\[
c(z) := \int_{\mathbb{R}^p} s(\beta,z) \,d\beta< \infty,
\]
then the posterior density is well defined (i.e., proper), intractable
and given by
\[
\pi(\beta\mid z) = \frac{s(\beta,z)}{c(z)} .
\]
For the time being, we assume that the posterior is indeed
proper.\vadjust{\goodbreak}

Let $U$ and $V$ be independent random variables such that $U \sim
\mathrm{N}(0,1)$ and $V \sim\operatorname{Exp}(1)$. Also, define $\theta=
\theta(r) = \frac{1-2r}{r(1-r)}$ and $\tau^2 = \tau^2(r) =
\frac{2}{r(1-r)}$. Routine calculations show that the random variable
$\theta V + \tau\sqrt{V} U$ has the asymmetric Laplace distribution
with parameter $r$. \citet{kozukoba2011} exploit this
representation to construct a DA algorithm as follows. For
$i=1,\ldots,m$, let $(Z_i,Y_i)$ be independent pairs such that
$Z_i|Y_i=y_i \sim\mathrm{N}(x_i^T \beta+ \theta y_i, y_i \tau^2)$ and,
marginally, $Y_i \sim\operatorname{Exp}(1)$. Then $Z_i - x_i^T \beta$ has
the asymmetric Laplace distribution with parameter $r$, as in the
original model. Combining this model with the flat prior on $\beta$
yields the augmented posterior density defined as
\[
\pi(\beta, y \mid z) = \frac{1}{c'(z)} \Biggl[ \prod_{i=1}^m
\frac{1}{\sqrt{2 \pi\tau^2 y_i}} \exp\biggl\{ - \frac{1}{2 \tau^2
y_i} ( z_i - x_i^T \beta- \theta y_i )^2 \biggr\}
e^{-y_i} I_{\mathbb{R}_+}(y_i) \Biggr] ,
\]
where $c'(z) = r^m (1-r)^m c(z)$. Of course, $\int_{\mathbb{R}_+^m}
\pi(\beta, y \mid z) \,dy = \pi(\beta\mid z)$. This leads to a DA
algorithm based on the joint density $\pi(\beta, y \mid z)$, which is
viable because, as we now explain, simulation from $\pi(\beta\mid y,
z)$ and $\pi(y \mid\beta, z)$ is straightforward.

As usual, define $X$ to be the $m \times p$ matrix whose $i$th row is
the vector~$x_i^T$. We assume throughout that $m \ge p$ and that $X$
has full column rank, $p$. Also, let $D$ denote an $m \times m$
diagonal matrix whose $i$th diagonal element is~$( \tau^2 y_i
)^{-1}$. A~straightforward calculation shows that
\[
\beta\mid y, z \sim\mathrm{N}_p ( \mu, \Sigma) ,
\]
where $\Sigma= \Sigma(y,z) = ( X^T D X )^{-1}$, and, letting
$l$ denote an $m \times1$ vector of ones,
\[
\mu= \mu(y,z) = ( X^T D X )^{-1} \biggl( X^T D z -
\frac{\theta}{\tau^2} X^T l \biggr) .
\]
Also, it's clear from the form of $\pi(\beta, y \mid z)$ that, given
$(\beta,z)$, the $y_i$'s are independent, and $y_i$ has density
proportional to
%
%e11 #&#
%
\begin{equation}
\label{eqymarginal}
\frac{1}{\sqrt{y_i}} \exp\biggl\{ -\frac{1}{2 \tau^2} \biggl[ y_i
(2\tau^2 + \theta^2) + \frac{( z_i - x_i^T \beta)^2}{y_i }
\biggr] \biggr\} I_{\mathbb{R}_+}(y_i) .
\end{equation}
This is\vspace*{1pt} the density of the reciprocal of an inverse Gaussian random
variable with parameters $2 + \theta^2/\tau^2$ and
$\sqrt{2\tau^2+\theta^2}/|z_i-x_i^T\beta|$. Thus, one iteration of
the DA algorithm requires one draw from a $p$-variate normal
distribution, and $m$ independent inverse Gaussian draws. Note that
in this example \mbox{$\X= \mathbb{R}^p$} and $\Y= \mathbb{R}_+^m$, so both
spaces have uncountably many points.

From this point forward, we restrict ourselves to the special case
where $r=1/2$, that is, to median regression. The proof of the
following result, which is fairly nontrivial, is provided in
Appendix \ref{appprop1}.
%
%pr1 #&#
%
\begin{proposition}
\label{propbqr}
If $r=1/2$ and $X$ has full column rank, then the joint density upon
which Kozumi and Kobayashi's DA algorithm is based satisfies~%
\eqref{eqct}. Thus, the corresponding Markov operator is trace
class.
\end{proposition}
%
%re2 #&#
%
\begin{remark}
Proposition \ref{propbqr} implies that, if $r=1/2$ and $X$ has full
column rank, then the posterior is proper, that is, $c(z)<\infty$.
First, by construction, the function $s(\beta,z)$ is an invariant
density for the DA Markov chain, whether it is integrable
(in~$\beta$) or not. Now, the fact that the DA Markov operator is trace
class implies that the DA Markov chain is geometrically ergodic,
which in turn implies that the chain is positive recurrent. Hence,
the chain cannot admit a nonintegrable invariant density
[\citet{meyntwee1993}, Chapter 10], so $s(\beta,z)$ must be
integrable, that is, the posterior must be proper.
\end{remark}

We now construct a sandwich algorithm for this problem. Let $G$ be
the multiplicative group, $\mathbb{R}_+$. Given $y \in\Y=
\mathbb{R}_+^m$ and $g \in\mathbb{R}_+$, define $g y$ to be scalar
multiplication of each element in $y$ by $g$, that is, $gy =
(gy_1,gy_2,\ldots,\allowbreak gy_m)$. Clearly, $e y = y$ and $(g_1 g_2)y = g_1
(g_2 y)$, so the compatibility conditions described in
Section \ref{secgroup} are satisfied. It is also easy to see that the
Lebesgue measure on $\Y$ is relatively invariant with multiplier
$\chi(g) = g^m$. When $r=1/2$, $\pi(y \mid z)$ is proportional to
\begin{eqnarray*}
&&\frac{e^{-\sum_{i=1}^m y_i}}{| X^T D X |^{1/2}}
\exp
\biggl\{ -\frac{1}{2} z^T D^{1/2} [ I - D^{1/2} X
(X^T D X)^{-1} X^T D^{1/2} ] D^{1/2} z \biggr\}\\
&&\qquad{}\times\prod_{i=1}^m y_i^{-{1/2}} I_{\mathbb{R}_+}(y_i) .
\end{eqnarray*}
Therefore, in this case, the density \eqref{eqgroupden} takes the
form
\begin{eqnarray*}
&&\frac{\pi(gy \mid z) g^m}{m(y)} \omega_l(dg)\\
&&\qquad\propto g^{({m-p-2})/{2}} e^{-g \sum_{i=1}^m y_i} \\
&&\qquad\quad{}\times\exp\biggl\{
-\frac{1}{2g} z^T D^{1/2} [ I - D^{1/2} X (X^T D
X)^{-1} X^T D^{1/2} ] D^{1/2} z \biggr\} \,dg.
\end{eqnarray*}
So at the middle step of the three-step procedure for simulating the
sandwich chain, we draw a $g$ from the density above and move from
$y=(y_1,y_2,\ldots,y_m)$ to $(gy_1,gy_2,\ldots,gy_m)$, which is a random
point on the ray that emanates from the origin and passes through the
point $y$. If $m$ happens to equal $p+1$, then this density has the
same form as \eqref{eqymarginal}, so we can draw from it using the
inverse Gaussian distribution as described earlier. Otherwise, we can
employ a simple rejection sampler based on inverse Gaussian and/or
gamma candidates. In either case, making one draw from this density
is relatively inexpensive.

Recall that $\pi(\beta\mid y, z)$ is a normal density. It's easy to
see that, if $g \ne1$, then $\pi(\beta\mid gy, z)$ is a different
normal density, which implies that \eqref{eqcanthold} does not hold.
Therefore, Theorems \ref{thmcompare} and \ref{thmatleastone} are
applicable and they imply that the ordered eigenvalues of the sandwich
chain are all less than or equal to the corresponding eigenvalues of
the DA chain, and at least one is strictly smaller. As far as we
know, this sandwich algorithm has never been implemented in practice.

%apA #&#
%
\begin{appendix}\label{app}

%seB #&#
\section{\texorpdfstring{Proof of Lemma \lowercase{\protect\ref{lemlambda1}}}{Proof of Lemma 1}}
\label{applemma1}

Fix $i \in\mathbb{N}$. Since $f(x,y) = f_X(x) f_Y(y)
\sum_{j=0}^\infty\beta_j g_j(x) h_j(y)$, we have
\[
(P_X h_i)(x) = \int_\Y h_i(y) \Biggl( \sum_{j=0}^\infty\beta_j
g_j(x) h_j(y) \Biggr) f_Y(y) \nu(dy) = \beta_i g_i(x) .
\]
A similar calculation shows that $P_Y g_i = \beta_i h_i$. Now, fix
$h \in L^2_0(f_Y)$. Because $\{h_i\}_{i=1}^\infty$ forms an
orthonormal basis for $L^2_0(f_Y)$, we have $h = \sum_{i=1}^\infty a_i
h_i$. Thus,
\[
\norm{P_X h} = \Biggl\| \sum_{i=1}^\infty a_i ( P_X h_i )
\Biggr\| = \Biggl\| \sum_{i=1}^\infty a_i \beta_i g_i \Biggr\| =
\sqrt{\sum_{i=1}^\infty a^2_i \beta^2_i} \le\beta_1 \norm{h} ,
\]
and we have equality if $h(y)=h_1(y)$. Hence, $\norm{P_X}_{L^2_0(f_Y)
\rightarrow L^2_0(f_X)} = \beta_1$. An analogous argument shows
that $\norm{P_Y}_{L^2_0(f_X) \rightarrow L^2_0(f_Y)} = \beta_1$. Now,
for each \mbox{$i \in\mathbb{N}$}, we have
\[
K g_i = P_X P_Y g_i = \beta_i P_X h_i = \beta^2_i g_i .
\]
But\vspace*{1pt} $\{g_i\}_{i=1}^\infty$ form an orthonormal basis of $L^2_0(f_X)$,
which proves that $K$ has eigenvalues $\{\beta_i^2\}_{i=1}^\infty$.

%seC #&#
\section{\texorpdfstring{Proof of Proposition \lowercase{\protect\ref{propbqr}}}{Proof of Proposition 1}}
\label{appprop1}

Here we show that the joint density underlying
Kozumi and Kobayashi's (\citeyear{kozukoba2011})
DA algorithm for median regression satisfies \eqref{eqct}. That is,
we will show that
\[
\int_{\mathbb{R}_+^m} \int_{\mathbb{R}^p} \pi(\beta\mid y, z)
\pi(y \mid\beta, z) \,d\beta \,dy < \infty.
\]
\begin{pf*}{Proof of Proposition \ref{propbqr}}
First,
\[
\pi(y \mid\beta, z) = c \frac{e^{-{a
y_{\cdot}}/{2}}}{\sqrt{\hat{y}}} \exp\Biggl\{
\frac{\sqrt{a}}{\tau} \sum_{i=1}^m |z_i-x_i^T\beta| - \frac{(z - X
\beta)^T D (z - X \beta)}{2} \Biggr\},
\]
where $a = (2 \tau^2 + \theta^2)/\tau^2$, $y_{\cdot} = \sum_{i=1}^m
y_i$, $\hat{y} = \prod_{i=1}^m y_i$, and $c$ is a constant (that does
not involve $y$ or $\beta$). Now let
\[
\mathcal{W} = \bigl\{ w \in\mathbb{R}^m \dvtx w_i \in\{-1,1\} \mbox{ for }
i=1,2,\ldots,m \bigr\} .
\]
For any $\beta\in\mathbb{R}^p$ and any $\sigma>0$, we have
\[
\exp\Biggl\{ \sigma\sum_{i=1}^m | z_i - x_i^T \beta|
\Biggr\} \le\exp\Biggl\{ \sigma\sum_{i=1}^m |z_i| \Biggr\} \sum_{w \in
\mathcal{W}} \exp\{ \sigma w^T X \beta\} .
\]
Thus, it suffices to show that, for every $w \in\mathcal{W}$,
\[
\int_{\mathbb{R}_+^m} \frac{e^{-{a
y_{\cdot}}/{2}}}{\sqrt{\hat{y}}} \biggl[ \int_{\mathbb{R}^p} \exp
\biggl\{ \frac{\sqrt{a}}{\tau} w^T X \beta- \frac{(z - X \beta)^T D
(z - X \beta)}{2} \biggr\} \pi(\beta\mid y, z) \,d\beta\biggr]
\,dy
\]
is finite. We start by analyzing the inner integral. First, recall
that $\pi(\beta\mid y, z)$ is a multivariate normal density with mean
$\mu= ( X^T D X )^{-1} X^T D z$ and variance $\Sigma= (
X^T D X )^{-1}$. Now,
\[
(z - X \beta)^T D (z - X \beta) = z^T D z + (\beta- \mu)^T
\Sigma^{-1} (\beta- \mu) - \mu^T \Sigma^{-1} \mu.
\]
Therefore, the integrand (of the inner integral) can be rewritten as
\begin{eqnarray*}
&&\exp\biggl\{ - \frac{1}{2} ( z^T D z - \mu^T \Sigma^{-1} \mu
) \biggr\} \exp\biggl\{ \frac{\sqrt{a}}{\tau} w^T X \beta
\biggr\} \frac{|2 X^T D X|^{1/2}}{(2
\pi)^{p/2}2^{p/2}} \\
&&\qquad{} \times\exp\biggl\{ -
\frac{1}{2} (\beta- \mu)^T 2 \Sigma^{-1} (\beta- \mu) \biggr\},
\end{eqnarray*}
so the inner integral can be expressed as
%
%e12 #&#
%
\begin{equation}
\label{eqeq1}\quad
\exp\biggl\{ - \frac{1}{2} ( z^T D z - \mu^T \Sigma^{-1} \mu
) \biggr\} \frac{1}{2^{p/2}} \int_{\mathbb{R}^p} \exp
\biggl\{ \frac{\sqrt{a}}{\tau} w^T X \beta\biggr\}
\tilde{\pi}(\beta\mid y, z) \,d\beta,
\end{equation}
where $\tilde{\pi}(\beta\mid y, z)$ is a multivariate normal density
with mean $\mu$ and variance~$\Sigma/2$. But the integral in
\eqref{eqeq1} is just the moment generating function of $\beta$
evaluated at the point $\sqrt{a} w^T X/\tau$. Hence, \eqref{eqeq1}
is equal to
\[
2^{-{p/2}} \exp\biggl\{ - \frac{1}{2} ( z^T D z - \mu^T
\Sigma^{-1} \mu) + \frac{\sqrt{a}}{\tau} ( w^T X \mu
) +
\frac{a}{4\tau^2} ( w^T X \Sigma X^T w ) \biggr\}.
\]
Now, straightforward manipulation yields
\[
z^T D z - \mu^T \Sigma^{-1} \mu= z^T D^{1/2} \bigl( I -
D^{1/2} X ( X^T D X )^{-1} X^T D^{1/2} \bigr)^2
D^{1/2} z \ge0 .
\]
It follows that $e^{- ( z^T D z - \mu^T \Sigma^{-1}
\mu)/2} \le1$. A similar calculation reveals that $w^T X
\Sigma\times\allowbreak
X^T w \le w^T D^{-1} w = \tau^2 y_{\cdot}$. Hence, \eqref{eqeq1} is
bounded above by
\[
2^{-{p/2}} \exp\biggl\{ \frac{\sqrt{a}}{\tau} ( w^T X
(
X^T D X )^{-1} X^T D z ) + \frac{a y_{\cdot}}{4} \biggr\}.
\]
Thus, it only remains to show that, for any $w \in\mathcal{W}$,\vspace*{2pt}
\[
\int_{\mathbb{R}_+^m} \frac{1}{\sqrt{\hat{y}}} \exp\biggl\{
-\frac{a
y_{\cdot}}{4} + \frac{\sqrt{a}}{\tau} ( w^T X ( X^T D X
)^{-1} X^T D z ) \biggr\} \,dy < \infty.\vspace*{2pt}
\]
We will prove this by demonstrating that $w^T X ( X^T D X
)^{-1} X^T D z$ is uniformly bounded in $y$.

It follows from the general matrix result established in
Appendix \ref{appmatrix} that, for each $i \in\{1,2,\ldots,m\}$ and
all $(y_1,y_2,\ldots,y_m) \in\mathbb{R}_+^m$,\vspace*{2pt}
\[
x_i^T \biggl( x_i x_i^T + \sum_{j \in\{1,2,\ldots,m\}, j \ne i}
\frac{y_i}{y_j} x_j x_j^T \biggr)^{-2} x_i \le C_i(X) ,\vspace*{2pt}
\]
where $C_i(X)$ is a finite constant. Thus,\vspace*{2pt}
\begin{eqnarray*}
\| ( X^T D X )^{-1} X^T D z \|_2 & = & \Biggl\|
\sum_{i=1}^m ( X^T D X )^{-1} \frac{x_i z_i}{\tau^2 y_i}
\Biggr\|_2 \\[2pt]
& \le &\sum_{i=1}^m \biggl\| ( X^T D X )^{-1}
\frac{x_i z_i}{\tau^2 y_i} \biggr\|_2 \\[2pt]
& = &\sum_{i=1}^m \Biggl\|
\Biggl( \sum_{j=1}^m \frac{x_j x_j^T}{\tau^2 y_j} \Biggr)^{-1}
\frac{x_i z_i}{\tau^2 y_i} \Biggr\|_2 \\[2pt]
& = &\sum_{i=1}^m |z_i|
\biggl\| \biggl(x_i x_i^T + \sum_{j \in\{1,2,\ldots,m\}, j \ne i}
\frac{y_i}{y_j} x_j x_j^T \biggr)^{-1} x_i \biggr\|_2 \\[2pt]
& = &
\sum_{i=1}^m |z_i| \sqrt{x_i^T \biggl( x_i x_i^T + \sum_{j \in
\{1,2,\ldots,m\}, j \ne i} \frac{y_i}{y_j} x_j x_j^T \biggr)^{-2}
x_i} \\[2pt]
& \le &\sum_{i=1}^m |z_i| C_i(X) .\vspace*{2pt}
\end{eqnarray*}
Hence,\vspace*{2pt}
\begin{eqnarray*}
| w^T X ( X^T D X )^{-1} X^T D z | & = &\| w^T X
( X^T D X )^{-1} X^T D z \|_2 \\ & \le &\| w^T X
\|_2 \| ( X^T D X )^{-1} X^T D z \|_2\vspace*{2pt}
\end{eqnarray*}
is uniformly bounded in $y$. This completes the proof.
\end{pf*}

%seD #&#
\section{A matrix result}
\label{appmatrix}

Fix $x_1,x_2,\ldots,x_n \in\mathbb{R}^p$ where $n$ and $p$ are
arbitrary positive integers. Now define
\begin{eqnarray*}
&&C_{p,n} (x_1; x_2,\ldots,x_n) \\
&&\qquad = \cases{
\displaystyle \sup_{a \in\mathbb{R}_+} x_1^T ( x_1 x_1^T + a I_p
)^{-2} x_1, &\quad if $n=1$,\vspace*{2pt}\cr
\displaystyle \sup_{a \in\mathbb{R}_+^n} x_1^T \Biggl( x_1 x_1^T +
\sum_{i=2}^n a_i x_i x_i^T + a_1 I_p \Biggr)^{-2} x_1, &\quad
if $n \geq2$.}
\end{eqnarray*}
%
%le3 #&#
%
\begin{lemma}
\label{lemmatrixresult} $C_{p,n}(x_1;x_2,\ldots,x_n) < \infty$.
\end{lemma}
\begin{pf}
We use induction on $p$. Note that when $p=1$, we have
\[
C_{1,n}(x_1;x_2,\ldots,x_n) = \sup_{a \in\mathbb{R}_+^n}
\frac{x_1^2}{(x_1^2 + \sum_{i=2}^n a_i x_i^2 + a_1 )^2} =
\cases{
0, &\quad $x_1 = 0$,\vspace*{2pt}\cr
\dfrac{1}{x_1^2}, &\quad $x_1 \ne0$,}
\]
which is finite in either case. Thus, the result is true for $p=1$.

Now assume that for any $n \in\mathbb{N}$ and any $x_1,\ldots,x_n \in
\mathbb{R}^{p-1}$,
\[
C_{p-1,n}(x_1;x_2,\ldots,x_n) < \infty.
\]
We will complete the argument by showing that, for any $n \in
\mathbb{N}$ and any $x_1,\ldots,x_n \in\mathbb{R}^p$,
$C_{p,n}(x_1;x_2,\ldots,x_n) < \infty$. The result is obviously true
when $x_1=0$. Suppose that $x_1 \ne0$, and let $P$ be an orthogonal
matrix such that $P x_1 = \norm{x_1}_2 e_1$, where $e_1 = (1, 0, 0,
\ldots, 0)^T \in\mathbb{R}^p$. For $i=2,3,\ldots,n$, define $b_i = P
x_i$. Then we have
\begin{eqnarray*}
&&x_1^T \Biggl( x_1 x_1^T  + \sum_{i=2}^n a_i x_i x_i^T + a_1 I_p
\Biggr)^{-2} x_1 \\
&&\qquad = x_1^T \Biggl( P^T P x_1 x_1^T P^T P +
\sum_{i=2}^n a_i P^T P x_i x_i^T P^T P + a_1 P^T P \Biggr)^{-2} x_1
\\
&&\qquad = x_1^T \Biggl( P^T \Biggl( \norm{x_1}_2^2 e_1 e_1^T +
\sum_{i=2}^n a_i b_i b_i^T + a_1 I_p \Biggr) P \Biggr)^{-2} x_1 \\
&&\qquad =
x_1^T P^{-1} \Biggl( \norm{x_1}_2^2 e_1 e_1^T + \sum_{i=2}^n a_i b_i
b_i^T + a_1 I_p \Biggr)^{-1} (P^T)^{-1} \\
&&\qquad\quad{}\times
P^{-1} \Biggl( \norm{x_1}_2^2 e_1 e_1^T + \sum_{i=2}^n
a_i b_i b_i^T + a_1 I_p \Biggr)^{-1} (P^T)^{-1} x_1 \\
&&\qquad =
\norm{x_1}_2^2 e_1^T \Biggl( \norm{x_1}_2^2 e_1 e_1^T + \sum_{i=2}^n
a_i b_i b_i^T + a_1 I_p \Biggr)^{-2} e_1 .
\end{eqnarray*}
Now let $A = \{ i \in\{2,\ldots,n\} \dvtx b_i^T e_1 = 0 \}$ and
let $B = \{2,\ldots,n\} \setminus A$, that is, $B = \{ i \in
\{2,\ldots,n\} \dvtx b_i^T e_1 \ne0 \}$. If $i \in A$, then there
exists a $v_i \in\mathbb{R}^{p-1}$ such that
\[
b_i = \left[
\matrix{
0 \cr
v_i}
\right]
\]
and, if $i \in B$, then there exists a nonzero real number $u_i$ and
$v_i \in\mathbb{R}^{p-1}$ such that
\[
b_i = \left[
\matrix{
u_i \cr
u_i v_i}
\right] .
\]
Thus, we have
\begin{eqnarray*}
&&
x_1^T  \Biggl( x_1 x_1^T + \sum_{i=2}^n a_i x_i x_i^T + a_1 I_p
\Biggr)^{-2} x_1 \\
&&\qquad = \norm{x_1}_2^2 e_1^T \Biggl( \norm{x_1}_2^2
e_1 e_1^T + \sum_{i=2}^n a_i b_i b_i^T + a_1 I_p \Biggr)^{-2} e_1
\\
&&\qquad = \norm{x_1}_2^2 e_1^T \biggl(
\left[
\matrix{
\norm{x_1}_2^2 & 0^T \cr
0 & 0}
\right]
+ \sum_{i \in A} a_i
\left[
\matrix{
0 & 0^T \cr
0 & v_i v_i^T}
\right]
\\
&&\qquad\quad\hspace*{59pt}{} + \sum_{i \in B} a_i u_i^2
\left[
\matrix{
1 & v_i^T \cr
v_i & v_i v_i^T}
\right]
+ a_1 I_p \biggr)^{-2} e_1
\\
&&\qquad = \norm{x_1}_2^2 e_1^T
\left[
\matrix{
u & v^T \cr
v & W
}
\right]^{-2} e_1 ,
\end{eqnarray*}
where $u := \norm{x_1}_2^2 + \sum_{i \in B} a_i u_i^2 + a_1$, $v :=
\sum_{i \in B} a_i u_i^2 v_i$ and
\[
W := \sum_{i \in A} a_i v_i v_i^T + \sum_{i \in B} a_i u_i^2 v_i v_i^T
+ a_1 I_{p-1} .
\]
If $B$ is empty, then $v$ is taken to be the zero vector in
$\mathbb{R}^{p-1}$. The formula for the inverse of a partitioned
matrix yields
\begin{eqnarray*}
\left[
\matrix{
u & v^T \cr
v & W}
\right]^{-1}
&=&
\frac{1}{u - v^T W^{-1} v}\\
&&{}\times\left[
\matrix{
1 & -v^T W^{-1} \cr
-W^{-1} v & (u - v^T W^{-1} v)
W^{-1} + W^{-1} v v^T W^{-1}}
\right] .
\end{eqnarray*}
It follows that
\[
e_1^T \left[
\matrix{
u & v^T \cr
v & W}
\right]^{-2} e_1 =
\frac{1+v^T W^{-2} v}{( u - v^T W^{-1} v )^2} .
\]
If $n=1$ or $B$ is empty, then
\[
C_{p,n} (x_1; x_2, \ldots, x_n) = \| x_1 \|_2^2 \sup_{a \in
\mathbb{R}_+^n} \frac{1}{( \| x_1 \|_2^2 + a_1 )^2} =
\frac{1}{\|x_1\|_2^2} < \infty,
\]
so the result holds. In the remainder of the proof, we assume that $n
\geq2$ and~$B$ is not empty.

Note that the matrix
\[
\left[
\matrix{
u - \norm{x_1}_2^2 & v^T \cr
v & W}
\right] = \sum_{i=2}^n a_i b_i b_i^T + a_1 I_p
\]
is positive definite, which implies that it's determinant is strictly
positive, that is,
\[
|W| ( u - \norm{x_1}_2^2 - v^T W^{-1} v ) > 0 .
\]
Since $W$ is also positive definite, $u - v^T W^{-1} v >
\norm{x_1}_2^2$. Moreover,
\[
v^T W^{-2} v = \norm{W^{-1} v}_2^2 = \biggl\|W^{-1} \biggl( \sum_{i
\in
B} a_i u_i^2 v_i \biggr) \biggr\|_2^2 \le\biggl[ \sum_{i \in B}
\| W^{-1} ( a_i u_i^2 v_i) \|_2 \biggr]^2 .
\]
Therefore,
\[
\frac{1+v^T W^{-2} v}{( u - v^T W^{-1} v )^2} \le\frac{1 +
[ \sum_{i \in B} \| W^{-1} ( a_i u_i^2 v_i)
\|_2 ]^2}{\norm{x_1}_2^4} .
\]
Putting all of this together yields
\[
x_1^T \Biggl( x_1 x_1^T + \sum_{i=2}^n a_i x_i x_i^T + a_1 I_p
\Biggr)^{-2} x_1 \le\frac{1 + [ \sum_{i \in B} \| W^{-1}
( a_i u_i^2 v_i) \|_2 ]^2}{\norm{x_1}_2^2} .
\]
Recall that $A \cup B = \{2,3,\ldots,n\}$. For fixed $i \in B$, let
$k_{i,1},k_{i,2},\ldots,k_{i,n-2}$ denote the $n-2$ elements of the set
$\{2,3,\ldots,n\} \setminus\{i\}$. Then we have
\begin{eqnarray*}
&&\| W^{-1}  ( a_i u_i^2 v_i) \|_2^2 \\
&&\qquad = v_i^T
\biggl( v_i v_i^T + \sum_{j \in A} \frac{a_j}{a_i u_i^2} v_j v_j^T +
\sum_{j \in B, j \ne i} \frac{a_j u_j^2}{a_i u_i^2} v_j v_j^T +
\frac{a_1}{a_i u_i^2} I_{p-1} \biggr)^{-2} v_i \\
&&\qquad \le
C_{p-1,n-1}(v_i;v_{k_{i,1}},v_{k_{i,2}},\ldots,v_{k_{i,n-2}}).
\end{eqnarray*}
Thus, using the induction hypothesis, we have
\begin{eqnarray*}
C_{p,n}(x_1;x_2,\ldots,x_n) & = & \sup_{a \in\mathbb{R}_+^n} x_1^T
\Biggl( x_1 x_1^T + \sum_{i=2}^n a_i x_i x_i^T + a_1 I_p \Biggr)^{-2}
x_1 \\
& \le &\frac{1 + [ \sum_{i \in B}
\sqrt{C_{p-1,n-1}(v_i;v_{k_{i,1}}, v_{k_{i,2}}, \ldots,
v_{k_{i,n-2}})} ]^2}{\norm{x_1}_2^2} ,
\end{eqnarray*}
which is finite. This completes the proof of the lemma.
\end{pf}
%
%re3 #&#
%
\begin{remark}
Note that if $x_1 x_1^T + \sum_{i=2}^n a_i x_i x_i^T$ is invertible
for every $(a_2,\ldots,\allowbreak a_n) \in\mathbb{R}_+^{n-1}$, then
\[
C_{p,n}(x_1;x_2,\ldots,x_n) = \sup_{(a_2,\ldots,a_n) \in
\mathbb{R}_+^{n-1}} x_1^T \Biggl( x_1 x_1^T + \sum_{i=2}^n a_i x_i
x_i^T \Biggr)^{-2} x_1 .
\]
\end{remark}
\end{appendix}

\section*{Acknowledgments}

The authors thank Jorge Rom\'{a}n and two anonymous
reviewers for helpful comments and suggestions.

%suskaldyti doi

% imsref loaded by lrinkeviciute, 2011-11-15 15:53:18
% imsref loaded by lrinkeviciute, 2011-11-15 16:04:44
%

\printaddresses

\end{document}